\documentclass[12pt]{article}
\usepackage{amssymb}
\usepackage{amsmath}
\usepackage{graphs}
\usepackage{amssymb}
\usepackage{amsthm}
\usepackage{url}
\usepackage[british]{babel}
\usepackage{bigstrut,array,float, fancyvrb,marvosym,verbatim,combgames}

\usepackage{psgo,array,url}

\theoremstyle{plain}
\newtheorem{definition}{Definition}

\newtheorem{theorem}[definition]{Theorem}

\newcommand{\lef}{\mathcal{L}}
\newcommand{\pre}{\mathcal{P}}
\newcommand{\n}{\mathcal{N}}
\newcommand{\ri}{\mathcal{R}}
\newcommand{\ti}{\mathcal{T}} 

\newcommand{\gfr}{G_F^{SR}}
\newcommand{\gfl}{G_F^{SL}}
\newcommand{\plus}{+_{\ell}}

\begin{document}
\title{Scoring Play Combinatorial Games}
\author{Fraser Stewart\\
\texttt{fraseridstewart@gmail.com}}
\date{}

\maketitle{}

\centerline{\bf Abstract}
In this paper we will discuss scoring play games.  We will give the basic definitions for scoring play games, and show that they form a well defined set, with clear and distinct outcome classes under these definitions.  We will also show that under the disjunctive sum these games form a monoid that is closed and partially ordered.  We also show that they form equivalence classes with a canonical form, and even though it is not unique, it is as good as a unique canonical form.

\section{Introduction}

Combinatorial game theory is the study of all two player perfect
information games.  It is the development of mathematical methods
that can be used to help understand the best strategies in playing
such games.  The full definition of a combinatorial game may be
found in \cite{LIP} or \cite{WW}.

There are two types of play in combinatorial game theory
\emph{normal} play, where the last player to move is the winner, and
\emph{mis\`ere} play, where the last player to move is the loser.

Mathematically we define combinatorial games recursively following
Conway \cite{ONAG}.

\begin{definition} The options of a game $G$ are the sets of all
games that a player can move to from $G$ and are denoted by:

$G^L=\{\hbox{All games that Left can move to from }G\}$

$G^R=\{\hbox{All games that Right can move to from }G\}$

A game $G$ is written as $\{G^L|G^R\}$ where $G^L$ and $G^R$ are the
options of Left and Right respectively.

\end{definition}

We abuse notation by letting $G^L$ and $G^R$ represent a set of
options and the specific options themselves.  The base case for the
recursion is the game $\{.|.\}$, i.e. $G^L=G^R=\emptyset$.  For the
rest of this paper we will represent the empty set with a ``$.$''.

Conway noticed that many games tend to break up into lots of smaller
ones, which are played all at once and independently. That is, a move
in one component will not affect any of the other components.  This
became the mathematical operator known as the disjunctive sum, and
is the operator which is most associated with combinatorial games.
It is defined as follows.

\begin{definition} The disjunctive sum of two games $G$ and $H$ is,

$$G+H=\{G^L+H,G+H^L|G^R+H,G+H^R\}.$$

\end{definition}

What this operator\footnote{We again abuse notation by allowing the
comma to mean set union, and the ``$+$'' means take the disjunctive
sum of the game with all games in the set, e.g. $G^L+H$ means take
the disjunctive sum of $H$ with all games in $G^L$.} says is that
when playing the game $G+H$, on your turn you may make a move on the
game $G$ or the game $H$ but not both.

When analysing games, we generally wish to know which \emph{outcome
class} a game is in.  The outcome class tells us who the winner of a
game will be when Left moves first and when Right moves first, but
only under alternate and optimal play.  The outcome classes are as
follows.

\begin{definition} We define the following:

\begin{itemize}
\item{$\mathcal{L}=\{G|\hbox{Left wins playing first or second in
}G$\}.}
\item{$\mathcal{R}=\{G|\hbox{Right wins playing first or second in
}G$\}.}
\item{$\mathcal{P}=\{G|\hbox{The second player to move wins in }G$\}.}
\item{$\mathcal{N}=\{G|\hbox{The first player to move wins in }G$\}.}
\end{itemize}

\end{definition}

A standard result is that these are the only four possible outcomes
under both normal and mis\`ere play \cite{LIP}, \cite{WW},
\cite{ONAG}.

Under normal and mis\`ere play there are several definitions which
allow us to determine when one move is better for a player than
another one.  They are as follows.

\begin{definition} We define the following:
\begin{itemize}
\item{$-G=\{-G^R|-G^L\}$.}
\item{For two games $G$ and $H$, $G=H$ if $G+X$ has the same outcome
as $H+X$ for all games $X$.}
\item{For two games $G$ and $H$, $G\geq H$ if Left wins $H+X$
implies that Left wins $G+X$ for all games $X$.}
\item{For two games $G$ and $H$, $G\leq H$ if Right wins $H+X$
implies that Right wins $G+X$ for all games $X$.}
\end{itemize}
\end{definition}

Under these definitions normal play games form a partially ordered
abelian group under the disjunctive sum, and are relatively easy to
analyse. On the other hand, mis\`ere play games do not form a group,
and while there are equivalence classes, they are
extremely difficult to analyse, by comparison to normal play games.

\section{Scoring Play Games}

Scoring play games are a third, and totally overlooked way of playing combinatorial games.  With scoring play we are no longer interested in who moves last, but who has accumulated the most points during play.  Examples of scoring play games are Go and Manacala, both of which are perhaps the oldest combinatorial games in existence.

In this paper, we will be looking at the structure of scoring play games under the disjunctive sum, since it is by far the most commonly used operator in combinatorial game theory.  Intuitively, we want all scoring play games to have the following four properties:

\begin{enumerate}
\item{The rules of the game clearly define what points are, and how players either gain or lose them.}
\item{When the game ends, the player with the most points wins.}
\item{For any two games $G$ and $H$, $a$ points in $G$ are equal to $a$ points in $H$, where $a\in \mathbb{R}$.}
\item{At any stage in a game $G$, if Left has $L$ points and Right has $R$ points then the score of $G$ is $L-R$, where $L,R\in \mathbb{R}$.}
\end{enumerate}

For example, in the game Go you get one point for each of your opponents stones that you capture, and for each piece of area you successfully take.  In Mancala you get one point for each bean you place in your Kala.  So when comparing these games, we would like one point in Mancala to be worth one point in Go.

Mathematically, scoring games are defined in the following way:

\begin{definition} A scoring play game $G=\{G^L|G^S|G^R\}$, where $G^L$ and $G^R$ are sets of games and $G^S\in\mathbb{R}$, the base case for the recursion is any game $G$ where $G^L=G^R=\emptyset$.\\

\noindent $G^L=\{\hbox{All games that Left can move to from } G\}$\\
\noindent $G^R=\{\hbox{All games that Right can move to from } G\}$,\\

and for all $G$ there is an $S=(P,Q)$ where $P$ and $Q$ are the number of points that Left and Right have on $G$ respectively.  Then $G^S=P-Q$, and for all $g^L\in G^L$, $g^R\in G^R$, there is a $p^L,p^R\in\mathbb{R}$ such that $g^{LS}=G^S+p^L$ and $g^{RS}=G^S+p^R$.

\end{definition}

A concept we will be using throughout this paper is the game tree of a game.  While it may be intuitively obvious to the reader, non-the-less, we feel it is important to define it mathematically.

\begin{definition}
The game tree of a scoring play game $G=\{G^L|G^S|G^R\}$ is a tree with a root node, and every node has children either on the Left or the Right, which are the Left and Right options of $G$ respectively.  All nodes are numbered, and are the scores of the game $G$ and all of its options.
\end{definition}

We also need to define a concept that we call the ``final score''.  This is something which hopefully the reader finds relatively intuitive.  When the game ends, which it will after a finite amount of time, the score is going to determine whether a player won, lost or tied.  

From a combinatorial game theory perspective we want to know ``what is the best that a player can do?''.  Left is trying to maximise the value of the score, while Right is trying to minimise it.  Since this is going to be the backbone of our theory it is important to get it right, and so we use the following definition.

\newpage

\begin{definition}  We define the following:

\begin{itemize}
\item{$\gfl$ is called the Left final score, and is the maximum score --when Left moves first on $G$-- at a terminal position on the game tree of $G$, if both Left and Right play perfectly.}
\item{$\gfr$ is called the Right final score, and is the minimum score --when Right moves first on $G$-- at a terminal position on the game tree of $G$, if both Left and Right play perfectly.}
\end{itemize}

\end{definition}

The reason we define it this way is because the terminal position can vary dramatically, depending on the rules of the game and the operator being used.  For instance under the long rule the game ends when a player cannot move on all components, but under the short rule the game ends when a player cannot move on any one component.  These two rules will clearly give different results when computing the final score of a game.

Since we want our definition to be as general as possible, i.e. cover every possibility, it makes sense to define the final score in this way.  For the purposes of this paper we will be using standard combinatorial game theory convention.  That is, a game ends when it is a player's turn and he has no options.

It is also important to note that we will only be considering finite games, i.e. for any game $G$ the game tree of $G$ has finite depth and finite width.  This means that $\gfl$ and $\gfr$ are always computable, and cannot be infinite or unbounded.  

There is also the case where a game may have a form of aggregate scoring.  For example players may play two games in sequence, and the winner would be the player who gets the most points over both games.  This gives scoring play games an additional dynamic, where in the event of a tie after two games, the winner may be determined by the player who managed to accumulate more points in one of the games.

However as far as this paper is concerned, we will not be considering games of this type.  We will only look at games where the winner is determined after one game ends.  Games with aggregate scoring would be an interesting area to look at for further research.

There are two conventions that we will be using throughout this paper.  The first is that in all examples given we will take the initial score of the game to be $0$, unless stated otherwise.  The second is that if for a game $G$, $G^L=G^R=\emptyset$, I will simply write $G$ as $G^S$, rather than $\{.|G^S|.\}$.   For example the game $G=\{\{.|0|.\}|1|\{.|2|.\}\}$, will be written as $\{0|1|2\}$.  The game $\{.|n|.\}$, will be written as $n$ and so on.  This is simply for convenience and ease of reading.

\subsection{An Example}\label{ch2ex}

Before we continue we will give an example of a scoring play game to demonstrate how to use the notation.  So consider the game Toad and Frogs from Winning Ways \cite{WW}, under scoring play.  The rules are as follows;

\begin{enumerate}
\item{The game is played on a horizontal grid.}
\item{Left moves Toads and Right moves Frogs.}
\item{Toads move from left to right and Frogs move from right to left.}
\item{Toads can only jump Frogs and Frogs can only jump Toads.}
\item{The player who jumps the most pieces wins.}
\end{enumerate}

\begin{figure}[htb]
\begin{center}
\begin{pspicture}(7,6)
\put(0,3){$TBF=$}\put(4,6){\cgtree{
  {TBF(0)}(+{BTF(0)}(|{FTB(-1)}({FBT(-1)}|)) |+{TFB(0)}({BFT(1)}(|{FBT(1)})|))
}}
\end{pspicture}
\end{center}
\caption{$TBF=\{\{.|0|\{-1|-1|.\}\}|0|\{\{.|1|1\}|0|.\}\}$}\label{tf1}
\end{figure}
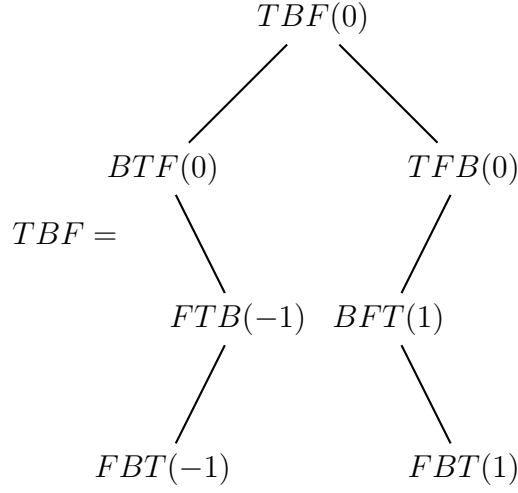

So consider the game $TBF$ as shown in figure \ref{tf1}, where $B$ represents a blank space, $T$ represents toads and $F$ represents frogs.  The numbers in brackets are the current score.

The game in figure \ref{tf1} has value $\{\{.|0|\{-1|-1|.\}\}|0|\{\{.|1|1\}|0|.\}\}$.  This game is in ``canonical form'', that is it neither has a dominated or reversible option.  For more details see section \ref{canfor}.

\subsection{Outcome Classes}

In combinatorial game theory we would like to know who wins under optimal play, e.g. if $G\in\lef$, then that means Left has a winning strategy moving first or second, if he plays his optimal strategy for both normal and mis\`ere play.  Under scoring scoring play the outcome classes are a little different, since in scoring play we allow ties, i.e. games where neither player wins.

Before we can define what the outcome classes precisely, we first need a new definition.  The definition we are about to give is very important for scoring play combinatorial game theory.  It, together with the definition of the final score, forms the core of our theory.

\begin{definition}
\item{$L_>=\{G|\gfl>0\}$, $L_<=\{G|\gfl<0\}$, $L_= =\{G|\gfl=0\}$.}
\item{$R_>=\{G|\gfr>0\}$, $R_<=\{G|\gfr<0\}$, $R_= =\{G|\gfr=0\}$.}
\item{$L_\geq=L_>\cup L_=$, $L_\leq = L_<\cup L_=$.}
\item{$R_\geq=R_>\cup R_=$, $L_\leq = R_<\cup R_=$.}
\end{definition}

Since we would like to classify every game by an outcome class it is also important that every game belongs to exactly one outcome class.  So we define the five outcome classes as follows.  

\begin{definition}
The outcome classes of scoring games are defined as follows:

\begin{itemize}
\item{$\lef=(L_>\cap R_>)\cup(L_>\cap R_=)\cup(L_=\cap R_>)$}
\item{$\ri=(L_<\cap R_<)\cup(L_<\cap R_=)\cup(L_=\cap R_<)$}
\item{$\n=L_>\cap R_<$}
\item{$\pre=L_<\cap R_>$}
\item{$\ti =L_=\cap R_=$}
\end{itemize}
\end{definition}

The reason that we chose the outcome classes in this way, is because if you have a game $G=\{1|0|0\}$, then it is more natural to say that belongs to the outcome $\lef$, since Right cannot win, but Left can if he moves first.  In this way we also keep the usual convention of calling a game $G\in \n$ a ``next player win'' and  a game $H\in \pre$ a ``previous player win''.

An interesting distinction is that while $\lef$ means the set of games where Left can win moving first or second in both normal and mis\`ere play, in scoring play, it means that if Left wins moving first he doesn't lose, and may win, moving second, and vice-versa.  Another distinction is the addition of the outcome class $\ti$, which of course does not exist in either normal or mis\`ere play, and means that the game ends in a tied score regardless of who moves first.

\begin{theorem}
Every game $G$ belongs to exactly one outcome class.
\end{theorem}

\begin{proof}  This is clear since ever game belongs to exactly one of $L_>$, $L_<$, $L_=$ and exactly one of $R_>$, $R_<$, $R_=$.  Therefore every game belongs to exactly one of the nine possible intersections of $L_>$, $L_<$, $L_=$ and $R_>$, $R_<$, $R_=$.  Since each outcome class is simply the union of one or more of these then each game can only be in exactly one outcome class.
\end{proof}

\subsection{The Disjunctive Sum}\label{disjunctive}

As we mentioned earlier, the disjunctive sum is by far the most commonly used operator in combinatorial game theory.  This is because many well known games such as Go naturally break up into the disjunctive sum of two or more components.  For scoring play the disjunctive sum needs to be defined a little differently, this is because in scoring play games when we combine them together we have to sum the games and the scores separately.  

For this reason we will be using two symbols $\plus$ and $+$.  The $\ell$ in the subscript stands for ``long rule'', this comes from \cite{ONAG}, and means that the game ends when a player cannot move on any component on his turn.  The ``short rule'' means that the game ends when a player cannot move on at least one component on his turn.  In this paper we will only be considering the disjunctive sum played with the long rule.

\begin{definition}  The disjunctive sum is defined as follows:
$$G\plus H=\{G^L\plus H,G\plus H^L|G^S+H^S|G^R\plus H,G\plus H^R\},$$
\noindent where $G^S+H^S$ is the normal addition of two real numbers.
\end{definition}

As with the disjunctive sum of normal and mis\`ere play games we abuse notation by making the comma mean set union, and $G^L\plus H$ means take the disjunctive sum of all $g^L\in G^L$ with $H$.

We would also like to know when one game is ``better'', than another one.  That is, given several options to play, which one is the best.  In normal play and mis\`ere play the definitions of ``$\geq$'' and ``$\leq$'' are relatively easy to define, since players either win or lose, however, for scoring play we have to take into account tied scores.  So for this reason we will re-define ``$\geq$'' and ``$\leq$''.

\begin{definition} We define the following:
\begin{itemize}
\item{$-G=\{-G^R|-G^S|-G^L\}$.}
\item{For any two games $G$ and $H$, $G=H$ if $G\plus X$ has the same outcome as $H\plus X$ for all games $X$.}
\item{For any two games $G$ and $H$, $G\geq H$ if $H\plus X\in O$ implies $G\plus X\in O$, where $O=L_\geq$, $R_\geq$, $L_>$ or $R_>$, for all games $X$.}
\item{For any two games $G$ and $H$, $G\leq H$ if $H\plus X\in O$ implies $G\plus X\in O$, where $O=L_\leq$, $R_\leq$, $L_<$ or $R_<$, for all games $X$.}
\item{$G\cong H$ means $G$ and $H$ have identical game trees.}
\item{$G\approx H$ means $G$ and $H$ have the same outcome.}
\end{itemize}
\end{definition}      

\begin{theorem}
$G\geq H$ if and only if $H\leq G$
\end{theorem}

\begin{proof}  First let $G\geq H$, and let $G\plus X\in O$ for some game $X$, where $O$ is one of $L_\leq$, $R_\leq$, $L_<$ or $R_<$.  This means that $H\plus X\not\in O'$, where $O'$ is one of $L_\geq$, $R_\geq$, $L_>$ or $R_>$, since if it was this would mean that $G\plus X\in O'$, since $G\geq H$, therefore $H\plus X\in O$, and hence $H\leq G$.

A completely identical argument can be used for $H\leq G$, and hence $G\geq H$ if and only if $H\leq G$ and the theorem is proven.
\end{proof}

\begin{theorem} Scoring play games are partially ordered under the
disjunctive sum.
\end{theorem}

\begin{proof}  To show that we have a partially ordered set we need
3 things.

\begin{enumerate}
\item{Transitivity: If $G\geq H$ and $H\geq J$ then $G\geq J$.}
\item{Reflexivity: For all games $G$, $G\geq G$.}
\item{Antisymmetry: If $G\geq H$ and $H\geq G$ then $G=H$.}
\end{enumerate}

1.  Let $G\geq H$ and $H\geq J$.  $G\geq H$ means that if $H\plus X\in O$ this implies $G\plus X\in O$, where $O=L_\geq$, $R_\geq$, $L_>$ or $R_>$, for all games $X$.  $H\geq J$, means that if $J\plus X\in O$ this implies that $H\plus X\in O$.  Since $G\geq H$, then this implies that $G\plus X\in O$, therefore $J\plus X\in O$ implies that $G\plus X\in O$ for all games $X$, and $G\geq J$.

2. Clearly $G\geq G$, since if $G\plus X\in O$ then $G\plus X\in O$, where  $O=L_\geq$, $R_\geq$, $L_>$ or $R_>$, for all games $X$.

3.  First let $G\geq H$ and $H\geq G$.  $G=H$ means that $G\plus X\approx H\plus X$ for all $X$.  So first let $G\plus X\in L_=$, then this implies that $H\plus X\in L_\geq$, since $H\geq G$.  However $H\plus X\in L_=$, since if $H\plus X\in L_>$, then this implies that $G\plus X\in L_>$, since $G\geq H$, therefore $G\plus X\in L_=$ if and only if $H\plus X\in L_=$.

An identical argument can be used for all remaining cases, therefore $G\plus X\approx H\plus X$ for all games $X$, i.e. $G=H$.
\end{proof}

\begin{theorem}\label{outcome}
For any three outcome classes $\mathcal{X}$, $\mathcal{Y}$ and $\mathcal{Z}$, there is a game $G\in\mathcal{X}$ and $H\in\mathcal{Y}$ such that $G\plus H\in\mathcal{Z}$.
\end{theorem}

\noindent\begin{proof}  Consider the games $G=\{\{\{d|c|e\}|b|.\}|a|.\}$ and $H=\{.|f|\{.|g|h\}\}$.  The final scores of $G$ are $\gfl=a$ and $\gfr=b$, and the final scores of $H$ are $H^{SL}_F=f$ and $H^{SR}_F=g$.  Now consider the game $G\plus H$ as shown in the figure.

\begin{figure}[htb]
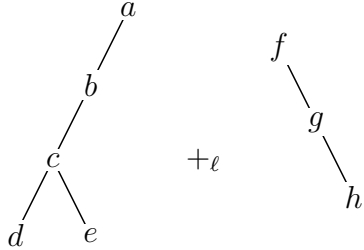

\begin{center}
\begin{graph}(4.5,3.5)

\roundnode{1}(0,0)\roundnode{2}(1,0)\roundnode{3}(0.5,1)\roundnode{4}(1,2)
\roundnode{5}(1.5,3)\roundnode{6}(3.5,2.5)\roundnode{7}(4,1.5)\roundnode{8}(4.5,0.5)

\edge{1}{3}\edge{2}{3}\edge{3}{4}\edge{4}{5}\edge{6}{7}\edge{7}{8}

\freetext(2.5,1){$\plus$}

\nodetext{1}{$d$}\nodetext{2}{$e$}\nodetext{3}{$c$}\nodetext{4}{$b$}
\nodetext{5}{$a$}\nodetext{6}{$f$}\nodetext{7}{$g$}\nodetext{8}{$h$}

\end{graph}
\end{center}
\caption{The game $G\plus H$, $G=\{\{\{d|c|e\}|b|.\}|a|.\}$ and $H=\{.|f|\{.|g|h\}\}$.}
\end{figure}

The final scores of $G\plus H$ are $(G\plus H)^{SL}_F=e+g$ or $d+h$ and $(G\plus H)^{SR}_F=e+h$.  Since $e$, $d$ and $h$ can take any value we can select them so that: $e+g$, $d+h$ and $e+h>0$ and $G\plus H\in \lef$; $e+g$, $d+h$ and $e+h<0$ and $G\plus H\in \ri$; $e+g$, $d+h>0$ and $e+h<0$ and $G\plus H\in \n$; $e+g$, $d+h<0$ and $e+h>0$ and $G\plus H\in \pre$ or finally $e+g=d+h=e+h=0$ and $G\plus H\in \ti$.

Since the outcomes of $G$ and $H$ depend on the values of $a$, $b$, $f$ and $g$, we can select them so that $G$ and $H$ can be in any outcome class, and thus the theorem is proven.

\end{proof}

Under normal play combinatorial games form an abelian group under the disjunctive sum.  The identity that is used is the set $\pre$, that is if $I\in \pre$ then $G\plus I\approx G$ for all games $G$.  In this case the entire set $\pre$ has a single unique representative, the game $\{.|.\}$.  This of course also means that $G=H$ if and only if $G\plus (-H)\in \pre$.

Under mis\`ere play, the identity set contains only one element, which is the same game $\{.|.\}$.  That is, if $G\not\cong \{.|.\}$, then $G\neq \{.|.\}$.  This was proven by Paul Ottaway \cite{O}.  This of course means that there is no easy or equivalent method for determining if two games are equivalent under mis\`ere play. 

For scoring play games, we have an equivalent theorem.  That is our identity set contains only one element, namely the game $\{.|0|.\}$, which we will call $0$.  It should be clear that $0\plus G\approx G$ for all games $G$, and so $0$ is the identity.

\begin{theorem}\label{identity}
For any game $G$, if $G\not\cong 0$ then $G\neq 0$.
\end{theorem}

\begin{proof}  The proof of this is very simple, first let $G^L\neq\emptyset$, since the case $G^R\neq\emptyset$ will follow by symmetry.  Next let $P=\{.|a|b\}$, and note that $P^{SL}_F=a$, since Left has no move on $P$.  So let $a>0$, if $G=0$ then this means that $(G\plus P)^{SL}_F\approx P$.  However, since $G$ is a combinatorial game we know from the definition that $G$ has both finite depth, and finite width.  So we can choose $b<0$ such that $|b|$ is greater than any score on the game tree of $G$.

Therefore when Left moves first on $G\plus P$ he must move to the game $G^L\plus P$.  Right will respond by moving to $G^L\plus b$, since $(G\plus P)^{SL}_F<0$ by choice of $b$.  This implies that $G\plus P\not\approx P$, and $G\neq 0$.  Hence the theorem is proven.
\end{proof}

What is interesting is that unlike mis\`ere games, some scoring games do have an inverse, namely the set of games $\{.|n|.\}$, where $n$ is a real number.  It should be clear that these are the only games which are invertible under scoring play, and any other non-trivial game cannot be inverted.

\section{Canonical Forms}\label{canfor}

Canonical forms are important, because if we can show that these games can be split up into equivalence classes with a unique representative for each class, then it makes these games much easier to analyze and compare.  We don't have to consider each game individually, but only the equivalence class to which it belongs.

\begin{theorem}\label{eq}
There exist two games $G$ and $H$ such that $G\not\cong H$ and
$G=H$.
\end{theorem}

\begin{proof} Consider the following games $G$ and $H$

\vspace{0.1in}

\begin{figure}[htb]
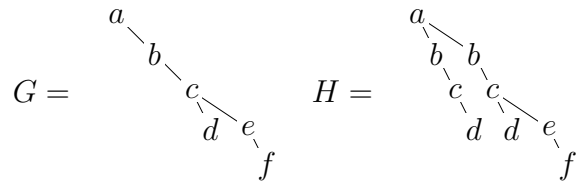

\begin{center}
\begin{graph}(7,2)

\graphnodesize{0.1} \graphlinewidth{0.01}

\roundnode{1}(1,2)\roundnode{2}(1.5,1.5)\roundnode{3}(2,1)\roundnode{4}(2.25,0.5)
\roundnode{5}(2.75,0.5)\roundnode{6}(3,0)\roundnode{7}(5,2)\roundnode{8}(5.25,1.5)
\roundnode{9}(5.75,1.5)\roundnode{10}(5.5,1.0)\roundnode{11}(6,1)\roundnode{12}(5.75,0.5)
\roundnode{13}(6.25,0.5)\roundnode{14}(6.75,0.5)\roundnode{15}(7,0)

\edge{1}{2}\edge{2}{3}\edge{3}{4}\edge{3}{5}
\edge{5}{6}\edge{7}{8}\edge{7}{9}\edge{8}{10}
\edge{10}{12}\edge{9}{11}\edge{11}{13}\edge{11}{14}\edge{14}{15}

\nodetext{1}{$a$}\nodetext{2}{$b$}\nodetext{3}{$c$}\nodetext{4}{$d$}
\nodetext{5}{$e$}\nodetext{6}{$f$}\nodetext{7}{$a$}\nodetext{8}{$b$}
\nodetext{9}{$b$}\nodetext{10}{$c$}\nodetext{11}{$c$}\nodetext{12}{$d$}
\nodetext{13}{$d$}\nodetext{14}{$e$}\nodetext{15}{$f$}

\freetext(0,1){$G=$}\freetext(4,1){$H=$}

\end{graph}
\end{center}
\caption{Two games $G$ and $H$, where $G\not\cong H$, but $G=H$.}
\end{figure}

\noindent where $a,b,c,d,e,f\in \mathbb{R}$.

This example is a variant of a similar example used to prove the
same theorem for mis\`ere games in \cite{GAMO}.

For any two games $G$ and $H$, $G=H$ if $G\plus X\approx H\plus X$ for all games $X$.  The easiest way to prove this is to show that $G\geq H$ and $H\geq G$.  Right can do at least as well playing $H\plus X$ as he can playing $G\plus X$, by simply copying his strategy from $G\plus X$ and not playing the left hand string on $H$.  Right cannot do better on $H\plus X$ than he can on $G\plus X$, since the string on the left hand side of $H$ can be copied on $G\plus X$ by simply not moving to $e$. So therefore if $H\plus X\in O$ then this implies that $G\plus X\in O$ where $O=L_\geq$, $R_\geq$, $L_>$ or $R_>$, i.e. $G\geq H$.

Left can also do at least as well playing $H\plus X$ as he can playing $G\plus X$, since if Right can achieve a lower final score playing the left hand string on $H\plus X$, then he can also do so by choosing not to move to $e$ on $G\plus X$.  Similarly if Right copies his strategy from $G\plus X$ onto $H\plus X$ then their final scores will be the same.  So if $G\plus X\in O$ then this implies that $H\plus X\in O$ where $O=L_\geq$, $R_\geq$, $L_>$ or $R_>$, i.e. $H\geq G$.  So therefore $G=H$ and the proof is finished.
\end{proof}

For both normal and mis\'ere play games, the standard way to reduce a game to its canonical form is to two concepts.  These are called domination and reversibility, and are defined as follows.

\begin{definition}
Let $G=\{A,B,C,\dots|G^S|D,E,F,\dots\}$, if $A\geq B$ or $D\leq E$ we
say that $A$ dominates $B$ and $D$ dominates $E$.
\end{definition}

\begin{definition}
Let $G=\{A,B,C,\dots|G^S|D,E,F,\dots\}$, an option $A$ is reversible if
$A^R\leq G$.  An option $D$ is also reversible if $D^L\geq G$.
\end{definition}

\begin{theorem}\label{dom}
Let $G=\{A,B,C,\dots|G^S|D,E,F,\dots\}$, and let $A\geq B$, then\\
$G'=\{A,C,\dots|G^S|D,E,F,\dots\}=G$.  By symmetry if $D\leq E$ and\\
$G''=\{A,B,C,\dots|G^S|D,F,\dots\}$ then $G''=G$.
\end{theorem}

\begin{proof}  Let $G=\{A,B,C,\dots|G^S|D,E,F,\dots\}$ such that $A\geq B$, further let\\ $G'=\{A,C,\dots|G^S|D,E,F,\dots\}$.  First suppose that $G\plus X\in O$, where $O=L_\geq$, $R_\geq$, $L_>$ or $R_>$ if Left moves to $B\plus X$.  This implies that $G'\plus X\in O$, since $A\geq B$.  Hence if $G\plus X\in O$  this implies that $G'\plus X\in O$, and since the Right options of $G$ and $G'$, this implies that $G'\geq G$.

Next suppose that $G'\plus X\in O'$ where $O'=L_\leq$, $R_\leq$, $L_<$ or $R_<$.  This implies that $G\plus X\in O'$, since the only option in $G^L$ that is not in $G'^L$ is $B$ and $B\leq A$, therefore $G'\leq G$, and $G=G'$.  So this means that the option $B$ may be disregarded and the proof is finished.

\end{proof}

\begin{theorem}\label{rev1}
Let $G=\{A,B,C,\dots|G^S|D,E,F,\dots\}$, and let $A$ be reversible with
Left options of $A^R=\{W,X,Y,\dots\}$.  If\\
$G'=\{W,X,Y,\dots,B,C,\dots|G^S|D,E,F,\dots\}$, then $G=G'$.  By symmetry  if $D$ is reversible with Right options of $D^L=\{T, S, R,\dots\}$.  If \\ $G''=\{A,B,C,\dots|G^S|T,S,R,\dots,D,E,F,\dots\}$, then $G=G''$.  
\end{theorem}

\begin{proof} Let $G=\{A,B,C,\dots|G^S|D,E,F,\dots\}$, where the Left options of $A^R=\{W,X,Y,\dots\}$ and let $G'=\{W,X,Y,\dots,B,C,\dots|G^S|D,E,F,\dots\}$, further let $A^R\leq G$.  If $G\plus X\in O$, where $O=L_\geq$, $R_\geq$, $L_>$ or $R_>$, when Left does not move to $A$ on $G$, then clearly $G'\plus X$ is also in $O$, since all other options for Left on $G$ are available for Left on $G'$.  

So consider the case where $G\plus X\in O$ if Left moves to $A\plus X$, then this implies that $A^R\plus X$ must also be in $O$.  This means that $G'\plus X\in O$ because $A^{RL}\subset G'^L$, and since all other options on $G'$ are the same as $G$, then $A^R\plus X\in O$ implies that $G'\plus X\in O$.  Hence if $G\plus X\in O$  then this implies that $G'\plus X\in O$, for all games $X$, i.e. $G'\geq G$.

Next assume that $G\plus X\in O'$, where $O'=L_\leq$, $R_\leq$, $L_<$ or $R_<$, for all games $X$.  However $A^R\leq G$, i.e. $G\plus X\in O'$ implies that $A^R\plus X\in O'$, and since $A^{RL}\subset G'^L$, and all other options on $G'$ are identical to options on $G$, this means that $G\plus X\in O'$, implies that $G'\plus X\in O'$, for all games $X$, i.e. $G'\leq G$.  Therefore $G=G'$ and the theorem is proven.
\end{proof}

\begin{definition}
A vertex $v$ on the game tree of a game $G$ is called a termination vertex if there is a game $X$, such that $G\plus X$ ends if the players reach vertex $v$.
\end{definition}

The reason why we called it a termination vertex is because it is a place where a game could potentially end.  If both Left and Right have an option at a particular vertex, then under the disjunctive sum a game cannot end at that point.

\begin{definition}
We say that $G$ is equivalent to $H$, or $G\equiv H$, if the underlying game trees of $G$ and $H$ are identical, and all termination vertices have the same score.
\end{definition}

\begin{theorem}
If $G\equiv H$, then $G=H$.
\end{theorem}

\begin{proof}  To prove this, first let $G\equiv H$ and let $G\plus X\in O$, where $O=L_>, L_\geq, R_>$ or $R_\geq$.  Since the underlying game trees of $G$ and $H$ are identical, and all termination vertices have the same score, then Left can do at least as well on $H\plus X$ simply by copying his strategy from $G\plus X$.  If he does so, then he will arrive at the same termination vertex on $H\plus X$ as he did on $G\plus X$, and therefore the games will end with identical scores.

Therefore, if $G\plus X\in O$ then this implies that $H\plus X\in O$ i.e., $H\geq G$.  By a totally symmetrical argument we also have that $G\geq H$.  So, $G=H$ and the theorem is proven.
\end{proof}

Equivalence is a little stronger than equality, and a little weaker than saying two games are identical.  The reason we need it is because it is possible for two games, say $G=\{1|1|1\}$ and $H=\{1|0|1\}$, to be equal to each other, not identical and neither has a dominated or reversible option.

However we still want to use domination and reversibility to achieve a ``canonical form'', so we will say that non-termination vertices are not important in the sense of determining the winner.  So while it is not a ``true'' canonical form in the sense that it is not necessarily unique, it is still useful for studying games.

\begin{definition}
A game $G$ is in canonical form if it has no dominated or reversible options.
\end{definition}

\begin{theorem}\label{cong}
For any two games $G$ and $H$ if $G=H$, and both $G$ and $H$ are in
canonical form, then $G\equiv H$.
\end{theorem}

\begin{proof} Let $G$ and $H$ be two games such that $G=H$ and neither $G$, nor $H$ has a dominated or reversible option.  

So first let $H\plus X\in O$, where $O=L_<, R_<, L_\leq$ or $R_\leq$.  Since $G=H$, this implies that $G\plus X\in O$.  However, if Left moves to $G^L\plus X$ then $G^{LR}\plus X$ cannot be in $O$.  If it were, this would mean that $H\plus X\in O$ implies $G^{LR}\plus X\in O$, i.e. $G^{LR}\leq H$, and $G$ would have a reversible option.  Which means that $G^L\plus X^R\in O$.  

This implies that $H^L\plus X^R\not\in O'$, where $O'=L_>, R_>, L_\geq$ or $R_\geq$.  Since if it were, then $H$ would have a dominated option.  Therefore, $G^L\plus X^R\in O$ if and only if $H^L\plus X^R\in O$, i.e. for all $g^L\in G^L$ there is an $h^L\in H^L$ such that $g^L\leq h^L$, and for all $h^L\in H^L$ there is a $g^{L'}\in G^L$ such that $h^L\leq g^{L'}$.

So that means $g^L\leq h^L\leq g^{L'}$.  However, $g^L$ and $g^{L'}$ must be identical, otherwise $g^L$ is a dominated option.  So every Left option of $G$ is equivalent to a Left option of $H$, i.e. $G^L\subseteq H^L$, and by a symmetrical argument $H^L\subseteq G^L$.  Therefore, $H^L\equiv G^L$, and similarly $H^R\equiv G^R$.

Since all options of $G$ and $H$ are equivalent, we can conclude that the only differences between the game trees of $G$ and $H$ are on non-terminating vertices. Therefore, $H\equiv G$ and the proof is finished. 
\end{proof}

It is also important to note that a game may have more than one canonical form.  For example, consider the game $G=\{\{3|0|4\},\{3|1|4\}|0|.\}$.  This game has two canonical forms, namely $\{\{3|0|4\}|0|.\}$ and $\{\{3|1|4\}|0|.\}$.  However both of these games are equivalent, so either can be used as the canonical form and it will not affect the analysis of this game.

\begin{theorem}
Let $G\to G_1\to G_2\to\dots \to G_n$ represent a series of
reductions on a game $G$ to a game $G_n$, which is in canonical
form.  Further let $G\to G_1'\to G_2'\to\dots\to G_m'$ represent a
different series of reductions on $G$ to a game $G_m'$ which is also
in canonical form, then $G_n\equiv G_m'$
\end{theorem}

\begin{proof}  Since each reduction preserves equality, then $G_n=G_m'$ and they are both in canonical form.  By theorem \ref{cong} $G_n\equiv G_m'$, and so the theorem is proven.
\end{proof}

Finally, it is important to note that it is certainly possible to define a unique canonical form, for example we could simply set all non-terminating vertices on a game tree to $0$.  However we feel that it is more important to keep the original values as this gives a lot of information about the games.  

Consider the games $G=\{3|0|4\}$ and $H=\{3|10|4\}$.  In the game $G$ Left moves and gains 3 points, while Right moves and loses 4, but in the game $H$ Left moves and loses 7 points, while Right moves and gains 6.  If we set $H^S$ to zero then this information would be lost.  So for this reason, we feel the only ways you should reduce a game is using domination and reversibility.
                                                                                                                                                                                                                             
\section{Acknowledgements}

I would like to thank my supervisor Keith Edwards for all his thoughtful comments and suggestions, as well as proof reading this paper.  He also gave me the idea to use the sets $L_>$, $L_<$ and so on.

\end{document}